\documentclass[11pt]{amsart}
\usepackage{epsfig}
\usepackage{graphicx, subfigure}
\usepackage{amsmath, amssymb}
\newtheorem{thm}{Theorem}[section]

\theoremstyle{definition}
\newtheorem{defn}[thm]{Definition}
\theoremstyle{remark}

\numberwithin{equation}{section}

\newcommand{\eps}{\varepsilon}

\begin{document}


\title[Geometry of Weak Stability Boundaries]{Geometry of Weak Stability Boundaries}

\author{E.\ Belbruno$^\dag$}
\address{Courant Institute of Mathematical Sciences, New York University, New York, New York 10003}
\email{belbruno@cims.nyu.edu}

\author{M.\ Gidea$^\ddag$}
\address{School of Mathematics, Institute for Advanced Study, Princeton, NJ 08540, USA \\ and Department of Mathematics, Northeastern Illinois University,
           5500 N.\ St.\ Louis Avenue, Chicago, IL 60625, USA}
           \email{mgidea@neiu.edu}
\author{F.\ Topputo}
\address{Aerospace Engineering Department, Politecnico di Milano,
           Via La Masa, 34, 20156, Milan, Italy}
\email{topputo@aero.polimi.it}

\thanks{$^\dag$ Research of E.B. was partially supported by NASA/AISR grant NNX09AK61G Program of SMD}
\thanks{$^\ddag$ Research of M.G. was partially supported by NSF grant  DMS-0635607.}

\begin{abstract}
The  notion of a weak stability boundary has been successfully used to design low energy trajectories from the Earth to the Moon. The structure of this boundary has been investigated in a number of studies, where partial results have been obtained.  We propose a generalization of the weak stability boundary. We prove analytically that, in the context of the planar circular restricted three-body problem,  under certain conditions on the mass ratio of the primaries and on the energy, the weak stability boundary about the heavier primary coincides with a branch of the  global stable manifold of the Lyapunov orbit about one of the Lagrange points.
\end{abstract}


\maketitle
\section{Introduction}

We consider the planar circular restricted three-body problem for a small mass ratio of the primaries. We give a general definition of the weak stability boundary set in the region of the heavier primary. We consider the global stable   manifold of the Lyapunov orbit about  the Lagrange point located between the primaries.
We prove analytically that, under restrictions on the energy, the weak stability boundary coincides with the branch of the global stable manifold in the region of the heavier primary.

The  concept of WSB was  introduced in \cite{belbruno1987,belbruno1990} to design low energy transfers from Earth to Moon, and subsequently applied to the rescue of the Japanese mission Hiten in 1991.\footnote{The GRAIL mission of NASA, arriving at the Moon on January 1, 2012, is using the same transfer as Hiten \cite{Lemonick12}.} (See also \cite{belbruno2004}.)  A particular feature of the `WSB method' useful for applications is that it allows the capture of a spacecraft into an elliptic orbit about the Moon, with specified eccentricity  of the ellipse, and with specified true anomaly at the capture.

There has been considerable work devoted to understand  the concept of WSB from the point of view of dynamical systems, and to enhance its applicability (see, e.g.,  \cite{garcia2007,toputo2009,BGT2010,SousaT2011b}). A remarkable property of the WSB is that,  in the context of the planar circular restricted three-body problem, for some range of energies, and under some topological conditions on the hyperbolic invariant manifolds associated to the libration points, the weak stability boundary points coincide with the points on the stable manifolds satisfying some additional conditions. This has been observed numerically in \cite{garcia2007}, and argued geometrically in \cite{BGT2010}.

The classical definition of the WSB is  as follows:  for each radial
segment emanating from the Moon, we consider trajectories
that leave that segment at the periapsis of an osculating ellipse
whose semi-major axis is a part of the radial segment; a
trajectory is called weakly $n$-stable if it makes $n$ full turns around the
Moon without going around the Earth, and it has  negative Kepler energy when it returns to
the radial segment; if the  trajectory is weakly $(n-1)$-stable but fails to be weakly $n$-stable, it is called weakly $n$-unstable; the points that make the transition from the weakly $n$-stable regime to the weakly  $n$-unstable regime are by definition the points of the WSB of order $n$.

We note that WSB points lie on different Hamiltonian energy levels. Also, the WSB is not an invariant set for the Hamiltonian flow. We remark that, since the stability/instability criteria, as described above, are concerned with the behavior of  trajectories for finite time, they inherently introduce  `artifacts', i.e., points with very similar trajectories that are categorized differently with respect to these criteria. See  \cite{BGT2010,SousaT2011a}.

In the present note, we propose a more general definition of the WSB. We remove the condition that the infinitesimal mass leaves the radial segment at the periapsis of an osculating ellipse
whose semi-major axis is a part of the radial segment.   We remove the condition on negative Kepler energy at the return. We define a point on the radial segment as being weakly $n$-stable provided that it makes $n$ turns around the primary, such that the distance from the infinitesimal mass to the primary measured along the trajectory does not get bigger than some critical distance. Otherwise the point is redeemed as unstable. (Some of these ideas are also suggested in \cite{SousaT2011a}.) The main result of this paper  is that the WSB points, which make the transition from the weakly stable to the weakly unstable regime, are the points on the stable manifold of the Lyapunov orbit for the corresponding energy level.

The argument for the main result is analytical, relying on  topological arguments and estimates from \cite{conley1971,McGehee,llibre1985}. For this reason, we deal with the WSB set  about the heavier primary (unlike in the WSB original setting).

An interesting aspect of the WSB method is that it uses `local' information on the dynamics, namely the return of trajectories to a surface of section about one of the primaries, to infer some `global' information on the dynamics, namely the existence of trajectories that execute transfers from one primary to the other.

\section{Background}
\subsection{The planar circular restricted three-body problem} \label{subsec:PCRTBP}
We consider the planar circular restricted three-body problem (PCRTBP) with the mass ratio of the primaries sufficiently small.  The system consists of two  mass points $P_1,P_2$, called primaries, of masses $m_1>m_2>0$, respectively, that  move under mutual
Newtonian gravity on circular orbits about their barycenter, and a third point $P_3$, of infinitesimal mass, that moves in the same plane as the primaries under their gravitational influence, but without exerting any influence on them.
Let  $\mu = m_2/(m_1 + m_2)$ be the relative mass ratio of $m_2$. In the sequel, we will assume that $0<\mu<1$ is very small, which will be made precise later.

It is customary to study the motion of the infinitesimal mass in a co-rotating system of coordinates $(x,y)$ that rotates with the primaries. Relative to this system, $P_1$ is positioned at $(\mu,0)$ and  $P_2$ is positioned at $(-1+\mu,0)$.   After some rescaling, the  equations of motions are given by
\begin{equation} \label{eqn:eq1}
\ddot{x} -2\dot{y}  = \displaystyle\frac{\partial\omega}{\partial
x},\qquad \ddot{y} +2\dot{x}  =
 \displaystyle\frac{\partial\omega}{\partial y},
\end{equation}
where the effective potential $\omega$ is given by
\begin{equation} \label{eqn:eq2}
\omega(x,y) = \frac{1}{2}(x^2+y^2) +
    \frac{1-\mu}{r_1} + \frac{\mu}{r_2} + \frac{1}{2}\mu(1-\mu),
\end{equation}
 with $r_{1} =
((x-\mu)^2 + y^2)^{1/2}$, $r_{2} = ((x+1-\mu)^2 + y^2)^{1/2}$.

The equations of motion can be described by a Hamiltonian system given by the following Hamiltonian (energy function):
\begin{equation} \label{eqn:ham}
H(x,y,p_x,p_y) = \frac{1}{2}((p_x+y)^2+(p_y-x)^2)-\omega(x,y),
\end{equation}
where $\dot x=p_x+y$ and $\dot y=p_y-x$.

\begin{figure} \centering
\includegraphics*[width=0.4\textwidth, clip, keepaspectratio]
{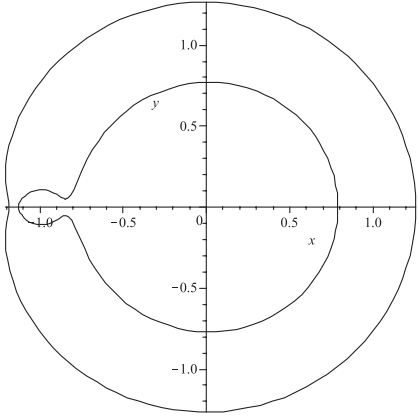}
    \caption[]{A Hill's region, $H\in(H(L_1),H(L_2))$.}
    \label{hill}
\end{figure}

For each fixed value $H$ of the Hamiltonian, the energy hypersurface $M_H$ is a  non-compact $3$-dimensional manifold in the $4$-dimensional phase space.
The projection of the energy hypersurface onto the configuration space $(x,y)$ is called a
Hill's region, and its boundary is a zero velocity curve. See Fig. \ref{hill}. Every
trajectory  is confined to the Hill's region corresponding to the
energy level of that trajectory.

The equilibrium points of the differential equations
\eqref{eqn:eq1} are given by the critical points of $\omega$. There
are five equilibrium points for this problem:  three of them, $L_1$, $L_2$ and $L_3$, are collinear with the
primaries (where $L_1$ is between $L_2$ and $L_3$), while the other two,
 $L_4, L_5$, form equilateral triangles with the
primaries. The distance from $L_1$ to $P_2$ is
given by the only positive solution $x_+$ to Euler's quintic equation
\begin{equation}\label{eqn:euler}x^5-(3-\mu)x^4+(3-2\mu)x^3-\mu x^2+2\mu x-\mu=0,\end{equation}
and so the distance from $L_1$ to $P_1$ is $1-x_+$.

The values   $H(L_i)$ of the Hamiltonian \eqref{eqn:ham} at the points $L_i$, $i=1,\ldots,5$, satisfy $H(L_5)=H(L_4)>H(L_3)>H(L_2)>H(L_1)$.
For $H<H(L_1)$, the Hill's  region has three components: two  bounded components, one about $P_1$ and the other  about $P_2$,
and a third component which is unbounded. For $H\in(H(L_1),H(L_2))$, the Hill's  region has two components, one bounded, which is topologically equivalent to the connected sum of the two bounded components from the case $H<H(L_1)$, and the other one unbounded (Fig. \ref{hill}).

The linearized stability of the equilibrium point $L_1$ is of saddle-center type, with the linearized equations possessing a pair of non-zero real eigenvalues $\pm \lambda$, and a pair of complex conjugate, purely imaginary eigenvalues $\pm i\nu$.  For each $H \gtrapprox H(L_1)$,  near the equilibrium point $L_1$ there exists a unique hyperbolic periodic orbit $\gamma_H$, referred as a Lyapunov orbit. This orbit has $2$-dimensional stable and unstable manifolds $W^s(\gamma_H)$, $W^u(\gamma_H)$, respectively, that are locally diffeomorphic  to $2$-dimensional cylinders. These manifolds have the following separatrix property: when restricted to a  compact  neighborhood $B_H(a,b)$ of $\gamma_H$ in the energy hypersurface $M_H$, of the type $B_H(a,b)=\{a\leq x\leq b\}$, with  $a<x_{L_1}<b$ sufficiently close to $x_{L_1}$,  each of the manifolds $W^s(\gamma_H), W^u(\gamma_H)$ separates $B_H(a,b)$ into two connected components.

\subsection{Conley's isolating block}\label{subsec:conley}

Let $\phi:M\times\mathbb{R}\to M$ be a $C^1$-flow on a $C^1$-differentiable manifold $M$. Given a  compact submanifold with boundary $B\subseteq M$, with $\dim(B)=\dim(M)$, we define
\begin{eqnarray*}
B^-=\{p\in \partial B \,|\, \exists \eps >0 \textrm { s.t. } \phi_{(0, \eps)}(p)\cap B=\emptyset \}, \\
B^+=\{p\in \partial B \,|\, \exists \eps >0 \textrm { s.t. } \phi_{(-\eps,0)}(p)\cap B=\emptyset  \},\\
B^0=\{p\in \partial B \,|\, \phi_t \textrm{ is tangent to }\partial B \textrm { at } p  \}.
\end{eqnarray*}

We obviously have $\partial B=B^0\cup B^-\cup B^+$. We call  $B^-$ the exit set and  $B^+$ the entry set of $B$.

An open set $V$ is called an isolating neighborhood for the flow if $\partial V$ contains no orbit of $\phi$.
 An invariant set $S$ for the flow $\phi$ is an isolated invariant set if there exists an isolating neighborhood $V$ for the flow such that $S$ is the maximal invariant set in $V$.
The compact submanifold $B$ is called an isolating block for the flow $\phi$ provided that:
\begin{itemize}\item[(i)] $B^-\cap B^+=B^0$, \item[(ii)] $B^0$ is a smooth
submanifold of $\partial B$ of codimension $1$, and, consequently, $B^-,B^+$ are submanifolds with common boundary $B^0$.\end{itemize}

The interior of an
isolating block is an isolating neighborhood and so determines an isolated invariant
set, possibly empty.

In the PCRTBP, Conley has constructed  an isolating block around $L_1$ that can be used to study the nearby dynamics.
Consider the part of the Hill's region which satisfies $a\leq x\leq b$, where $(a, b)$ contains the $x$-coordinate
$x_{L_1}$ of $L_1$. This set determines a ``dynamical channel'' which allows for the transit of trajectories between the $P_1$ and $P_2$ regions. The lift $B_H=B_H(a,b)$ of this set to  the energy hypersurface, where $a,b$ are chosen close to $x_{L_1}$,   is Conley's isolating block. Geometrically, this is a $3$-dimensional manifold with boundary $\partial B_H$ consisting of the set of points in the energy hypersurface that projects onto $x = a$ and $x=b$ in the configuration space.
It is diffeomorphic to the product of a line segment with a two sphere, $B_H \approx [a,b]\times S^2$, and its boundary
$\partial B_H$   is diffeomorphic to the union of two $2$-spheres, $\partial B_H= B_{H,a}\cup  B_{H,b}\approx(\{a\}\times S^2)\cup (\{b\}\times S^2)$.

The isolating block conditions in this case are that
every trajectory intersecting $\partial B$ tangentially must lie outside of $B_H$ both before and after the
intersection, that is,
if  $x(t) = a$  and $\dot x (t) = 0$ then $\ddot x(t) < 0$  and if $x(t) = b$, and $\dot x (t) = 0$ then $\ddot x(t) > 0$.
So we have
\begin{eqnarray*}B_H^0 &=& \{(x,y,\dot x,\dot y)\in \partial B_H\,|\,  x(t) = a \textrm{ or }  x(t) = b \textrm{ and } \dot x (t) = 0\}, \\
B_H^-&=&\{(x,y,\dot x,\dot y)\in \partial B_H\,|\,  x(t) = a \textrm{ and } \dot x (t) < 0, \textrm { or }  x(t) = b \textrm{ and } \dot x (t) > 0\}, \\
B_H^+&=&\{(x,y,\dot x,\dot y)\in \partial B_H\,|\,  x(t) = a \textrm{ and } \dot x (t) > 0, \textrm { or }  x(t) = b \textrm{ and } \dot x (t) < 0\}.
\end{eqnarray*}

For  each component of $\partial B_H$, the exit and entry sets determine  a pair of disjoint open $2$-dimensional topological disks, which we denote as follows:
$B_{H,a}^{-}$,  $B_{H,b}^{-}$ are  the exit sets of the boundary components $B_{H,a}$, $B_{H,b}$, respectively,
 and $B_{H,a}^{+}$, $B_{H,b}^{+}$  are the entry sets of the boundary components $B_{H,a}$, $B_{H,b}$, respectively. The complement in $B_{H,b} $ of $B^{-}_{H,b}\cup  B^{+}_{H,b}$ is the set $B_{H,b}^{0}=B_H^0\cap \{x=b\}$. A similar statement holds for  $B_{H,a}$.

\begin{figure} \centering
\includegraphics*[width=0.7\textwidth, clip, keepaspectratio]
{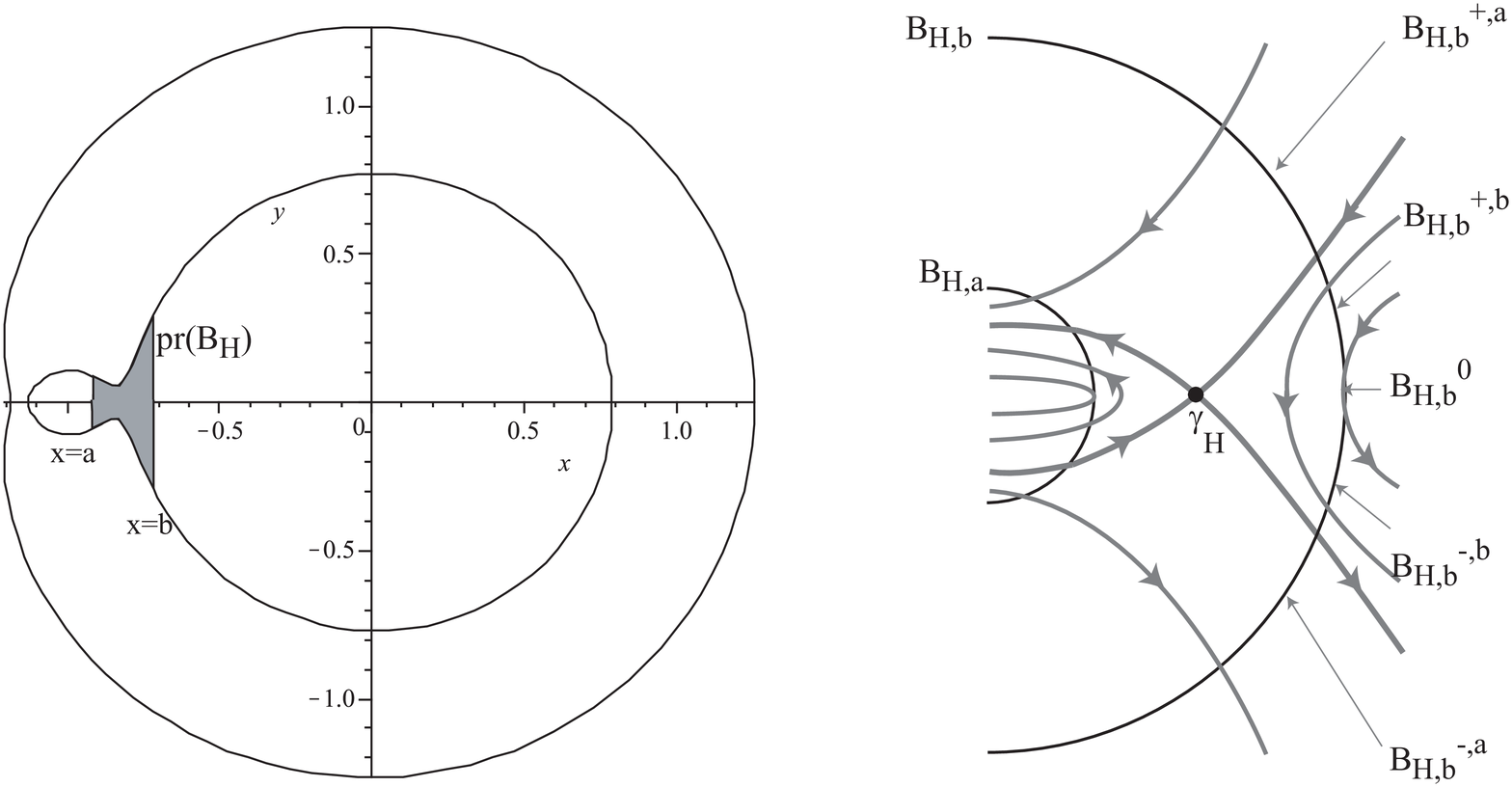}
    \caption[]{(a) Projection of Conley's isolating block onto configuration space. (b) Schematic representation of the dynamics across Conley's isolating block.}
    \label{conley}
\end{figure}

The exit and entry sets are further broken up into components with dynamical roles. The set $B_{H,b}^{+}$ is the union of three sets, a spherical cap $B_{H,b}^{+,a} $, corresponding to trajectories that enter the block $B_H$ through the entry part of $B_{H,b}$ and later leave the block through the exit part of $B_{H,a}$, a spherical zone $B_{H,b}^{+,b} $, corresponding to trajectories that enter the block $B_H$ through the entry part of $B_{H,b}$ and leave the block through the exit part of $B_{H,b}$, and a topological circle separating them, corresponding to the intersection of $W^s(\gamma_H)$ with $B_{H,b}$. Similarly, $B_{H,b}^{-}=B_{H,b}^{-,a} \cup B_{H,b}^{-,b} \cup (B_{H,b} \cap W^u(\gamma_H))$, where the notation is analogous to the above. There is a similar decomposition for the entry and exit set components of  $B_{H,a}$. See Fig. \ref{conley}.

Later in the paper, we will use the following fact, which is a consequence of the above discussion. There are three possible behaviors for trajectories that start from the $P_1$-region and  enter the isolating block:
\begin{itemize}
\item [(i)] Trajectories  enter the block through $B_{H,b}^{+,a} $, exit the block through $B_{H,a}^{-,b} $, and so they execute a transfer from the $P_1$-region to the $P_2$-region.
    \item [(ii)] Trajectories   enter the block through $B_{H,b}^{+,b}$, exit the block through $B_{H,b}^{-,b} $, and so they do not transfer to the $P_2$-region.
    \item[(iii)] Trajectories   enter the block through $B_{H,b}\cap W^s(\gamma_H)$ and are forward asymptotic to $\gamma_H$, and so they never leave the block.
\end{itemize}

For further details on this subsection, see \cite{conley1971}.

\subsection{Hyperbolic invariant manifolds}\label{subsec:invman}

The geometry of the hyperbolic invariant manifolds can be described analytically inside the $P_1$-region, for some range of energies and mass ratios, following some results from  \cite{McGehee,llibre1985}.

First, there exists an open set $O_1$ in the $(\mu,H)$-parameter plane, with $0<\mu\ll 1$ and $H \gtrapprox H(L_1)$ such  that, for $(\mu,H)\in O_1$, the following hold:
\begin{itemize}
\item [(i)] The energy hypersurface $M_{H}$ contains an invariant $2$-torus $\mathcal{T}_H$ separating   $P_1$  from $L_1$.
\item [(ii)] There exist $a<x_{L_1}<b$ such that the flow inside the isolating block $B_H=B_H(a,b)$ is conjugate to the linearized flow.
\item [(iii)] In the region $\mathcal{N}_H$ in $M_H$ bounded by $\mathcal{T}_H$ and $B_{H,b}$, the longitudinal angular coordinate $\theta$ is increasing along trajectories.
\end{itemize}

Second, for all $0<\mu\ll 1$ sufficiently small, the $(x,y)$-projections of the branches of $W^u(L_1),W^s(L_1)$ inside the $P_1$-region have the following properties:
\begin{itemize}
\item [(iv)] The distance $d$ to the zero velocity curve, and the angular coordinate $\theta$ satisfy the following estimates:
\begin{eqnarray}
\label{eqn:llibre1} d&=\mu^{1/3}\left (\frac{2}{3}N-3^{1/6} +M\cos t + o(1) \right ),\\
\label{eqn:llibre2} \theta&=-\pi+\mu^{1/3}\left (Nt+2M\sin t+o(1) \right ),
\end{eqnarray}
where $M,N$ are constants, the parameter $t$ means the physical time measured from a suitable origin, and  $o(1)\to 0$ when $\mu\to 0$ uniformly in $t$ as $t=O(\mu^{-1/3})$.   These expressions hold true outside $B_H$.
\item[(v)] There exists an open set $O_2\subseteq O_1$ in the $(\mu,H)$-parameter plane, with $0<\mu\ll 1$ and $H \gtrapprox H(L_1)$ such  that, for $(\mu,H)\in O_2$, the $(x,y)$-projections of the branches of $W^u(\gamma_H), W^s(\gamma_H)$ inside the $P_1$-region satisfy estimates similar to \eqref{eqn:llibre1} and \eqref{eqn:llibre2}. That is, these invariant manifolds turn around $P_1$ in the region $\mathcal{N}_H$ bounded by the torus $\mathcal{T}_H$ and the boundary component $B_{H,b}$ of the isolating block $B_H$.
    Moreover, there exists a sequence of mass ratios $\mu_k$ for which $W^u(\gamma_H)$ and $W^s(\gamma_H)$ have symmetric transverse intersections, provided $(\mu_k,H)\in O_2$.
\end{itemize}

The geometry of the hyperbolic invariant manifolds for the range of parameters considered above allows to extend the separatrix property of these manifolds from the local case, as described in Subsection \ref{subsec:PCRTBP}, to the global case. For as long as the stable and unstable manifolds do not intersect each other, the cuts of these manifolds with a surface of section are topological circles. If a point is inside the $i$-th cut $\Gamma_{\theta_0,i}^s(\gamma_H)$ made by the stable manifold $W^s(\gamma_H)$ with the surface of section $S_{\theta_0}$, which is assumed to be a topological circle,  then the forward trajectory of that point stays inside the cylinder bounded by $W^s(\gamma_H)$ in $M_H$ for $i$-turns and transfers from the $P_1$-region to the $P_2$-region afterwards. If a point in $S_{\theta_0}$ is outside the $i$-th cut $\Gamma_{\theta_0,i}^s(\gamma_H)$, then its forward trajectory stays inside the $P_1$-region for at least $(i+1)$-turns.
A similar statement holds for the cuts made by the unstable manifold and backwards trajectories.

If the stable and unstable manifolds intersect, say $\Gamma_{\theta_0,i}^s(\gamma_H)$ intersects $\Gamma_{\theta_0,j}^u(\gamma_H)$, then the intersection points are homoclinic points that make $(i+j)$-turns about $P_1$, and some future cuts of the invariant manifolds cease to be topological circles.  For example, $\Gamma_{\theta_0,i+j}^s(\gamma_H)$ is a finite union of open curve segments whose endpoints wind asymptotically towards $\Gamma_{\theta_0,i}^s(\gamma_H)$. Due to the asymptotic behavior of the endpoints, each of these open curves divides  $S_{\theta_0}$ into transfer and non-transfer orbits. Thus, the separatrix property extends to the case when the cuts of the hyperbolic invariant manifolds cease being topological circles.
See \cite {gideamasdemont2007,BGT2010}.

There are no analogues of the above analytical results for the $P_2$-region about the lighter mass.

\begin{figure} \centering
\includegraphics*[width=0.4\textwidth, clip, keepaspectratio]
{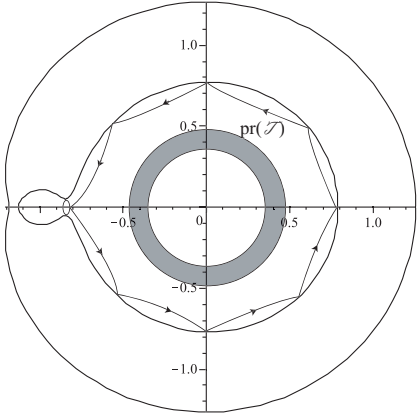}
    \caption[]{Projection of McGehee's separating torus onto configuration space, and trajectory near the zero velocity curve.}
    \label{mcgehee}
\end{figure}

\subsection{Equations of motion relative to polar coordinates}

We recall the relations between the motion of the infinitesimal mass $P_3$  relative relative to the barycentric
rotating coordinates $(x, y,\dot x, \dot y)$, relative to the polar coordinates $(r,\theta, \dot r,\dot \theta)$ about $P_1$, and relative to the classical orbital elements  $(a, e, \phi, \tau)$ about $P_1$.

The relation between barycentric and polar coordinate is  $r=((x-\mu)^2+y^2)^{1/2}$ and $\tan\theta=y/(x-\mu)$.

The orbital elements are characterized by the semi-major axis  $a$ of an ellipse with a focus at $P_1$,  the ellipse eccentricity $e\in[0,1)$, the argument of the periapsisis $\phi\in[0,2\pi]$, and the true anomaly $\tau\in [0,2\pi]$.
We have the following coordinate
transformations
\begin{equation} \label{eq:t1}
    \begin{split}
x & =   r \cos(\phi+\tau) + \mu, \\
y & =   r \sin(\phi+\tau),       \\
\dot x   &=  \dot r\cos(\phi+\tau) - r\dot\tau\sin(\phi +
\tau) + r \sin(\phi+\tau), \\ \dot y & =   \dot
r\sin(\phi+\tau) + r\dot\tau\cos(\phi + \tau) - r
\cos(\phi+\tau),
    \end{split}
\end{equation}
 and the following formulas
\begin{equation} \label{eq:t2}
\begin{split}
r & =  \frac{a(1 - e^2)}{1+ e \cos\tau}, \\
\dot r & =   \frac{a e (1-e^2)\dot\tau\sin\tau}{(1+e\cos \tau)^2},\\
\theta&=\phi+\tau,\\
\dot\theta & =\dot\tau=   \frac{\sqrt{a(1-e^2)(1-\mu)}}{r^2}.
    \end{split}
\end{equation}

The Hamiltonian function in polar coordinates is given by
\begin{equation}\label{eqn:hampol}H(r,p_r,\theta,p_\theta)= \frac{1}{2}(p_r^2+\frac{1}{r^2}\theta^2)-p_\theta+\mu r\cos\theta +\omega(r,p_r,\theta,p_\theta),\end{equation}
where the canonical momenta are given by
\[p_r=\dot r,\, p_\theta=r^2(\dot\theta+1).\]

Note that the conservation of energy implies that the initial position $(r,\theta)$ relative to $P_1$  and the initial radial velocity $\dot r$ uniquely determine a trajectory, up to a choice of a sign for $\dot\theta$. Suppose that we know the initial data $(r,\dot r,\theta)$ on a trajectory. Using \eqref{eq:t2}, the  eccentricity of the osculating ellipse to this trajectory at the initial point uniquely determines the trajectory, and hence its energy. This implicitly defines $\phi$ and $\tau$. Conversely, if we have a trajectory for which the initial angle coordinate $\theta$, the initial angular velocity $\dot r$, and the eccentricity of the osculating ellipse $e$ at the initial condition are fixed, then the energy level $H$ of the trajectory uniquely determines its initial value of $r$.

\section{Weak Stability Boundary}

We consider the system of polar coordinates $(r,\theta)$ about $P_1$ as above, and we let $H(r,\dot r,\theta,\dot \theta)$ be the Hamiltonian relative to this coordinate system. As discussed above, the energy is also uniquely determined by the $(r,\dot r,\theta, e)$-data, where $e$ is the eccentricity of the osculating ellipse at the initial point.
We consider a Poincar\'e section through $P_1$ that makes an angle $\theta_0$ with the $x$-axis, which is given by
\[S_{\theta_0}=\{(r,\dot r,\theta,\dot \theta)\,|\,
\theta=\theta_0,\, \dot\theta>0\}.\]

Let $l_{\theta_0}$ denote the radial segment obtained as the intersection of $S_{\theta_0}$ with the $(x,y)$-space.
Any trajectory that meets $S_{\theta_0}$ transversally is uniquely determined by the $(r,\dot r)$-coordinates of the intersection point,
as the $\theta$-coordinate equals $\theta_0$ in this section, and the $\dot
\theta$-coordinate can be solved uniquely from the energy condition
$H(r,\dot r,\theta,\dot \theta)=H$, provided  $\dot\theta>0$.

Consider a trajectory with the initial condition $z_0=z_0(r_0, \dot r_0,\theta_0,e_0)$ with initial position $r(0)=r_0$, $\theta(0)=\theta_0$, initial radial velocity $\dot r(0)=\dot r_0$, and $\dot\theta(0)>0$, for which the osculating ellipse at the initial point has eccentricity $e_0$. We   keep the values of $\dot r_0, \theta_0, e_0$ fixed and investigate the change of behavior of the trajectories   when $r_0$ changes.
Note that different initial values of $r_0$ yield different energies $H_0$.

Fix a  value of $\mu$ sufficiently small for which there exists an open range of energies $(H(L_1),H^*)$ with $(\mu, H)\in O_2$ for each $H\in (H(L_1),H^*)$, as in Subsection \ref{subsec:invman}. For this range of energies the estimates \eqref{eqn:llibre1} are valid.
Fix $a<x_{L_1}<b$ such that $B_H(a,b)$ is an isolating block for all $H\in(H(L_1),H^*)$.  Let $y_b$ be the supremum of the $y$-coordinates on the segment $x=b$ inside the Hill's regions for $H\in (H(L_1),H^*)$. Define  $\theta_1=\arctan({y_b}/(\mu-b))$. Let $D_1$ be the distance from $P_1$ to $x=a$, that is $D_1=\mu-a$.

Fix  $H\in(H(L_1),H^*)$ and consider the projection $\textrm{pr}_{(x,y)}(\mathcal{N}_H)$ of $\mathcal{N}_H$ onto the $(x,y)$-configuration plane. For  each angle coordinate $\theta\in[0,2\pi]$, there exists a well defined  interval  $(r_1(H,\theta), r_2(H,\theta))$ such that $(r,\theta)\in \textrm{pr}_{(x,y)}(\mathcal{N}_H)$ if and only if $r\in (r_1(H,\theta), r_2(H,\theta))$.
For each trajectory point $(r,\theta) \in \textrm{pr}_{(x,y)}(\mathcal{N}_H)$ there exists a set of admissible values   of the radial velocity $\dot r$ and of the eccentricity of the osculating ellipse $e$ corresponding to the trajectory at that point. When we let $H$ vary in $(H(L_1),H^*)$, then for each $\theta\in[0,2\pi]$, we obtain an open set of admissible values of $(\dot r,e)$ corresponding to all trajectories for all of these energy levels.

We fix an angle $\theta_0$ and a pair of admissible values $(\dot r_0,e_0)$. Since the energy $H$ is uniquely determined by the data $(r_0,\dot r_0,\theta_0,e_0)$, there exists an open set $\mathcal{R}(\dot r_0,\theta_0,e_0)\subseteq  (r_1(H,\theta), r_2(H,\theta))$ of $r_0$-values such that ${H_0}=H(r_0,\dot r_0,\theta_0,e_0)\in(H(L_1),H^*)$ provided $r_0\in \mathcal{R}(\dot r_0,\theta_0,e_0)$. In the next definition, we will consider trajectories with initial points $z_0$ lying on the radial segment $l_{\theta_0}$. We will restrict to values of $r_0$ in the set  $\mathcal{R}(\dot r_0,\theta_0,e_0)$.

\begin{defn}\label{defn:stable} We say that a forward trajectory with initial point $z_0=z_0(r_0,\theta_0)$ in $l_{\theta_0}$,    initial  radial velocity $\dot r_0$ and initial eccentricity of the osculating ellipse $e_0$, is weakly  $n$-stable provided that it turns $n$-times around  $P_1$, with all   intersections with $l_{\theta_0}$ being transverse,  and  such that the distance to $P_1$ is always less than  $D_1$. If the trajectory is weakly $(n-1)$-stable but fails to be weakly $n$-stable, we say that the trajectory is weakly $n$-unstable.
\end{defn}

The conditions on the parameters assumed for the  Definition \ref{defn:stable} are imposed in order to define the critical distance $D_1$ in a consistent way for the whole range of energy values $H\in(H(L_1),H^*)$. We recall that in the classical definition of the WSB, a trajectory is called $n$-stable  if it  turns $n$-times around  $P_1$, without turning around $P_2$; in that case    one can consider the distance from $P_1$ to $P_2$ as the critical distance.

We note that the transversality requirement in  Definition \ref{defn:stable}, on  the
intersections of the trajectory of the infinitesimal mass with $l_{\theta_0}$, implies that weak $n$-stability is an open condition, that is, if a trajectory starting at some  $z_0=(r_0,\dot r_0, \theta_0,e_0)$ is weakly $n$-stable,  then all
trajectory starting inside some domain of the type \[(r,\dot r, \theta, e)\in(r_0-\eps,r_0+\eps)\times (\dot r_0-\eps,\dot r_0+\eps)\times (\theta_0-\eps,\theta_0+\eps)\times (e_0-\eps,e_0+\eps)\] with $\eps>0$ sufficiently small, are also weakly $n$-stable.

Thus  we obtain the following  set of weakly stable points in the phase space
\[ {\mathcal{W}}_n=
\{z_0(r_0,\dot r_0,\theta_0,e_0)\,|\, z_0\textrm { is weakly $n$-stable relative to } l_{\theta_0},\, \theta_0\in[0,2\pi]\}.
\]

Due to the open conditions on the $n$-stable trajectories, the set ${\mathcal{W}}_n$ is an open set of points in the phase space.
If we fix the parameters $\dot r_0$, $\theta_0$ and $e_0$, then we obtain an open set ${\mathcal{W}}_n(\dot r_0,\theta_0,e_0)$ in $l_{\theta_0}$, which is a countable union of disjoint open intervals
\begin{equation}\label{eqn:r-star}{\mathcal{W}}_n(\dot r_0,\theta_0,e_0)=\bigcup_{k\geq 1}(r_{2k-1},r_{2k}).\end{equation}

The points of the type $r_{2k-1},r_{2k}$ at the ends  of these intervals   are weakly $n$-unstable.

\begin{defn}\label{defn:wsb}
The WSB of order $n$, denoted 
$\mathcal{W}^*_n$, is  the set  of all points $r^*(r_0,\dot r_0,\theta_0,e_0)$ that are at the boundary of the set of the weakly $n$-stable
points, i.e.,
\begin{displaymath}
\mathcal{W}^*_n =\partial {\mathcal{W}}_n.\end{displaymath}\end{defn}

We also  denote by $\mathcal{W}^*_n(\dot r_0,\theta_0,e_0)$  the set of WSB points on the radial segment $l_{\theta_0}$ of fixed parameters $\dot r_0$ and $e_0$. Thus, the WSB set $\mathcal{W}^*_n(\dot r_0,\theta_0,e_0)$ contains  the closure of the set of all points of the type $r_{2k-1},r_{2k}$,  which are the endpoints  of the intervals of weakly $n$-stable points within each radial segment $l_{\theta_0}$ as in \eqref{eqn:r-star}.

The main result of the paper says that, if we restrict to some angle range of $\theta_0$ outside the angle sector $[\pi-\theta_1,\pi+\theta_1]$, where $\theta_1$ is defined as above,  then the WSB set is completely determined by the stable manifolds of Lyapunov orbits. To state this result, we have to adopt a convention on how to count the number of cuts made by the stable manifold with a surface of section $S_{\theta_0}$. We label a cut made by the stable manifold $W^s(\gamma_{H_0})$ with $S_{\theta_0}$ as the $i$-th cut provided that the net change $\Delta \theta$ of the angle $\theta$ along all trajectories starting from $S_{\theta_0}$ and ending asymptotically at $\gamma_{H_0}$ satisfies $2i\pi\leq\Delta\theta<2(i+1)\pi$. Note that as long as  $\theta_0  \not\in[\pi-\theta_1,\pi+\theta_1]$ there is no ambiguity about the labeling of the cuts with the section $S_{\theta_0}$.

\begin{thm}\label{thm:main} Fix a pair of admissible values $(\dot r_0, e_0)$ as defined above. Assume  $\theta_0\in (-\pi+\theta_1,\pi-\theta_1)$, where $\theta_1$ is defined as above.
Then a point $z_0=z_0(r_0,\dot r_0,\theta_0,e_0)$, with $r_0\in\mathcal{R}(\dot r_0,\theta_0,e_0)$, is in $\mathcal{W}^*_n(\dot r_0,\theta_0,e_0)$ if and only if $z_0$ lies on the $(n-1)$-st cut $\Gamma^s_{\theta_0,n-1}(\gamma_{H_0})$ of the stable manifold $W^s(\gamma_{H_0})$ with the surface of section $\mathcal{S}_{\theta_0}$, where ${H_0}$ is the energy level corresponding to $z_0$.
\end{thm}

The restrictions  imposed on the parameters in Theorem \ref{thm:main}  are needed to apply the analytical arguments from Subsection  \ref{subsec:invman}.
It is nevertheless  shown in  \cite{BGT2010} that the WSB overlaps with some subset of  the stable manifold  of the Lyapunov orbit under much weaker conditions, provided that the hyperbolic invariant manifolds satisfy some topological condition (they turn around the primaries for a long enough time, without colliding with the primaries). Moreover, in \cite{BGT2010} a wider energy range is considered, in which case the WSB is identified with a subset of the union of the  stable manifolds of the Lyapunov orbits about  $L_1$ and about $L_2$.   The situation described by Theorem \ref{thm:main} is just a special case when the required topological conditions can be verified analytically.

Now we explain the relation between WSB and hyperbolic invariant manifolds   in a more concrete way. Assume that we fix some energy level ${H_0}\in(H(L_1),H^*)$. We generate the stable manifold of the Lyapunov orbit $\gamma_{H_0}$, and we count the successive cuts made by the stable manifold with some Poincar\'e surface of section $S_{\theta_0}$.  Let $z_0$ be a point on the $(n-1)$-st cut $\Gamma^s_{\theta_0,n-1}(\gamma_{H_0})$ of  $W^s(\gamma_{H_0})$ with  $\mathcal{S}_{\theta_0}$. Let $\dot r_0$ be the radial velocity at $z_0$, and $e_0$ the eccentricity of the osculating ellipse at $z_0$. Then the point $z_0$ is in the WSB set $\mathcal{W}^*_n(\dot r_0,\theta_0,e_0)$. Moreover, every WSB point can be obtained in this way.

\section{Proof of the main result}

Due to the angle restriction $\theta_0\in (-\pi+\theta_1,\pi-\theta_1)$, in Theorem \ref{thm:main}, we restrict to the following set of weakly $n$-stable points
\[\tilde{\mathcal{W}}_n= \{z_0(r_0,\dot r_0,\theta_0,e_0)\,|\, z_0\textrm { is weakly $n$-stable relative to } l_{\theta_0},\, \theta_0\in(-\pi+\theta_1,\pi-\theta_1)\}.\]

Since the $n$-stability is an open condition and the angle range $(-\pi+\theta_1,\pi-\theta_1)$ is also open, the set $\tilde{\mathcal{W}}_n$ is an open set in  the phase space.

We prove  that a point $z_0=z_0(r_0,\dot r_0,\theta_0,e_0)$ is in  $\tilde{\mathcal{W}}^*_n$ if and only if it is in the $(n-1)$-st cut $\Gamma^s_{\theta_0,n-1}(\gamma_{H_0})$ made by the stable manifold $W^s(\gamma_{H_0})$ with  $\mathcal{S}_{\theta_0}$, where ${H_0}$ is the energy level corresponding to $z_0$.
For this, we first  show that $z_0$ is a weakly  $n$-stable point on $l_{\theta_0}$ if and only if it is outside the domain in $S_{\theta_0}$  bounded by $\Gamma^s_{\theta_0,n-1}(\gamma_{H_0})$, and is  weakly $n$-unstable if and only if it is inside the domain in $S_{\theta_0}$  bounded by  $\Gamma^s_{\theta_0,n-1}(\gamma_{H_0})$.

First, we show that the points inside the cylinder bounded by the stable manifold are weakly unstable. Let  $z_0=z_0(r_0,\dot r_0,\theta_0,e_0)$ be a point in $l_{\theta_0}$. Then \eqref{eqn:hampol} gives the value ${H_0}$ of the energy of the trajectory with initial condition $z_0$. Assume that $z_0$ is inside the domain in $S_{\theta_0}$  bounded by $\Gamma^s_{\theta_0,n-1}(\gamma_{H_0})$. By the separatrix property from Subsection \ref{subsec:invman} the trajectory turns counterclockwise precisely $(n-1)$-times  inside the domain $\mathcal{N}_{H_0}$, while staying inside the region of the cylinder bounded by $W^s(\gamma_{H_0})$, enters the isolating block $B_{H_0}$ through the entry set region $B_{{H_0},b}^{+,a}$, crosses the block  and exits it through the exit set region $B_{{H_0},a}^{-,b}$. When the trajectory
leaves the block $B_{H_0}$, the distance from $P_1$ is bigger than $D_1$. Since the trajectory achieves a distance to $P_1$ bigger than the threshold value prior to completing an $n$-th turn around $P_1$, the trajectory is weakly $n$-unstable.

Second, we show that the points outside the cylinder bounded by the stable manifold are weakly stable. Assume that $z_0$ is outside  the domain in $S_{\theta_0}$  bounded by $\Gamma^s_{\theta_0,n-1}(\gamma_{H_0})$. By the separatrix property from Subsection \ref{subsec:invman} the trajectory will turn  counterclockwise inside the domain $\mathcal{N}_{H_0}$ and  will keep  staying outside the region of the cylinder bounded by $W^s(\gamma_{H_0})$ for at least $n$ turns. If the trajectory leaves the domain $\mathcal{N}_{H_0}$, it has to meet the block $B_{H_0}$ at $B^0_{H_0}$ or at $B_{{H_0},b}^{+,b}$. In the first case, the trajectory bounces back to the domain $\mathcal{N}_{H_0}$ and it continues its counterclockwise motion about $P_1$. In the second case, it cannot leave the block $B_{H_0}$ through $B_{{H_0},a}$, since only the points that are inside the cylinder bounded by $W^s(\gamma_{H_0})$ can do that; it cannot remain inside the block   $B_{H_0}$  for all future times since only the points on $W^s(\gamma_{H_0})$ have this property; hence, it has to leave $B_{H_0}$ through the exit set region $B_{{H_0},b}^{-,b}$, and to go back to the domain $\mathcal{N}_{H_0}$. The time spent by the trajectory inside the block $B_{H_0}$ does not affect the  count of turns about $P_1$.  Since $z_0$ is outside the cut $\Gamma^s_{\theta_0,n-1}(\gamma_{H_0})$, the trajectory cannot leave the $P_1$-region after only $(n-1)$-turns, so it  turns around $P_1$ for at least $n$-turns. Thus the trajectory is weakly $n$-stable.

Now, we prove the statement of the main theorem.

First, assume that $z_0$ is in $\Gamma^s_{\theta_0,n-1}(\gamma_{H_0})$. Its forward trajectory turns $(n-1)$-times around $P_1$ and then approaches asymptotically  $\gamma_{H_0}$. Thus the trajectory is weakly $n$-unstable.  To show that $z_0$ is an WSB point it is sufficient to prove that there exists a sequence $(z_0^k)_{k\geq 1}$ with $z_0^k\in  \tilde{\mathcal{W}}_n$ and $z_0^k\to z_0$ as $k\to\infty$.  Take a small $4$-dimensional open ball $\mathcal {U}$ around $z_0$ in the phase space. Let $T>0$ be  such that the time-$T$ map $\phi_T$ of the Hamiltonian flow  takes $z_0$ to a point in $B_{H_0,b}$, where $H_0=H(z_0)$.
The image $\phi_T(\mathcal{U})$ of $\mathcal{U}$ by $\phi_T$ is a $4$-dimensional open topological ball about  $\phi_{T}(z_0)$. We intersect $\phi_T(\mathcal{U})$ with the $4$-dimensional submanifold  with boundary $\bigcup_{H\in(H(L_1),H^*)}B_{H}$. The ball $\phi_T(\mathcal{U})$  has non-empty intersection with  $\bigcup_{H\in(H(L_1),H^*)}B^{+,b}_{H,b}$. These intersection points yield weakly $n$-stable trajectories.
Thus $\phi_T(\mathcal{U})\cap \bigcup_{H\in(H(L_1),H^*)}B_{H}$  contains a $4$-dimensional open, topological ball $\mathcal{V}$, which contains  $\phi_{T}(z_0)$ on its boundary, consisting of points that correspond to weakly $n$-stable trajectories, i.e., those trajectories that return to the $P_1$ region for at least one extra turn about $P_1$.

Now consider the  set $\phi_{-T}(\mathcal{V})$. This  is a $4$-dimensional open, topological ball in $\mathcal{U}$ that contains $z_0$ on its boundary. There exist $\theta'$ arbitrarily close to $\theta_0$ such that the intersection $\phi_{-T}(\mathcal{V})\cap S_{\theta'}$ is  a non-empty open set. All points $z'\in \phi_{-T}(\mathcal{V})\cap {S}_{\theta'}$ are weakly $n$-stable points. We note that these points may not lie on $l_{\theta_0}$, nor on the  same energy level as $z_0$; they can also have the eccentricity of the osculating ellipse  different from $e_0$.
Thus, arbitrarily near $z_0$ one can always find weakly $n$-stable points,  and since $z_0$ itself is weakly $n$-unstable, it follows that  $z_0\in \partial \tilde{\mathcal{W}}_n=\tilde{\mathcal{W}}^*_n$.

Second, assume that $z_0\in \tilde{\mathcal{W}}^*_n(\dot r_0,\theta_0, e_0)$. Then there exists a sequence of points $(z_0^k)_{k\geq 1}$ on $l(\theta_0)$ such that $z_0^k$ is weakly $n$-stable and $z_0^k\to z_0$ as $k\to \infty$. From the above, we know that the  weakly $n$-stable points are those inside the cylinder bounded by the stable manifold.  Thus, there exists a corresponding sequence of stable manifold cuts $\Gamma^{s}_{\theta_0,n-1}(\gamma_{H_k})$ where $H_k=H(z_0^k)$, such that $z_0^k$ is inside the region in $S_{\theta_0}$ bounded by $\Gamma^{s}_{\theta_0,n-1}(\gamma_{H_k})$. Since $z_k^0\to z_0$ it follows that $H(z_0^k)\to H(z_0)=H_0$ as $k\to\infty$. The stable manifold cuts also depend continuously on the energy, so $\Gamma^{s}_{\theta_0,n-1}(\gamma_{H_k})$ approaches $\Gamma^{s}_{\theta_0,n-1}(\gamma_{{H_0}})$ as  $k\to\infty$. Hence $z_0\in \Gamma^{s}_{\theta_0,n-1}(\gamma_{{H_0}})$.

Through double inclusion, we conclude  that \[z_0\in \tilde{\mathcal{W}}^*_n(\dot r_0,\theta_0, e_0)\textrm{ if and only if }  z_0\in \Gamma^{s}_{\theta_0,n-1}(\gamma_{H_0}).\]

\section{Concluding remarks}

We compare the invariant manifold method with the WSB method. The invariant manifold method is based on identifying geometric objects that serve as building blocks that organize the global dynamics: equilibrium points, periodic orbits, and their stable and unstable invariant manifolds, if they exist. The WSB method is a local method for deciding whether the trajectories about one of the primaries exhibit some  kind of stability in terms of the return to a surface of section. The conclusion of this paper, corroborated with the results in \cite{garcia2007,BGT2010}, is that in simple models the two methods  overlap for a substantial range of parameters.

One can think of some other kinds of indicators that mark the passage between the weakly $n$-stable and the weakly $n$-unstable regimes. One such a possible indicator is the continuity of the Poincar\'e return map. The $n$-th return map to $S_{\theta_0}$ is continuous at all weakly $n$-stable points. At the WSB points, the return map exhibits essential discontinuities of infinite type. Thus, the set of points where the return map fails to be continuous contains the WSB points.

\end{document}